\documentclass[12pt]{amsart}

\usepackage{times,epsf,amssymb,amsmath,hyperref, rlepsf, color,qtree,tikz, wrapfig,mathtools}

\usepackage{tikz-qtree}
\usetikzlibrary{trees} 

\begin{document}

\newtheorem{thm}{Theorem}[section]
\newtheorem{lem}[thm]{Lemma}
\newtheorem{cor}[thm]{Corollary}
\newtheorem{conj}[thm]{Conjecture}
\newtheorem{question}[thm]{Question}

\theoremstyle{definition}
\newtheorem{defn}[thm]{\bf{Definition}}

\theoremstyle{remark}
\newtheorem{rmk}[thm]{Remark}

\def\square{\hfill${\vcenter{\vbox{\hrule height.4pt \hbox{\vrule width.4pt height7pt \kern7pt \vrule width.4pt} \hrule height.4pt}}}$}

\newenvironment{pf}{{\it Proof:}\quad}{\square \vskip 12pt}

\title[Minimal Surfaces in Hyperbolic $3$-manifolds]{Minimal Surfaces in Hyperbolic $3$-manifolds}
\author{Baris Coskunuzer}
\address{UT Dallas, Dept. Math. Sciences, Richardson, TX 75080}
\email{coskunuz@utdallas.edu}

\thanks{The author is partially supported by Simons Collaboration Grant, and Royal Society Newton Mobility Grant.}


\newcommand{\Si}{S^2_{\infty}(\mathbf{H}^3)}
\newcommand{\PI}{\partial_{\infty}}
\newcommand{\BH}{\mathbf{H}}
\newcommand{\BR}{\mathbb R}
\newcommand{\BC}{\mathbb C}
\newcommand{\BZ}{\mathbb Z}
\newcommand{\BN}{\mathbb N}

\newcommand{\N}{\mathbf{N}}

\newcommand{\e}{\epsilon}
\newcommand{\wh}{\widehat}
\newcommand{\wt}{\widetilde}

\newcommand{\D}{\mathcal{D}}
\newcommand{\I}{\mathcal{I}}
\newcommand{\F}{\mathcal{F}}
\newcommand{\A}{\mathcal{A}}
\newcommand{\U}{\mathcal{U}}
\newcommand{\B}{\mathbf{B}}
\newcommand{\C}{\mathcal{C}}
\newcommand{\X}{\mathcal{X}}
\newcommand{\Y}{\mathcal{Y}}
\newcommand{\E}{\mathcal{E}}
\newcommand{\s}{\mathcal{S}}
\newcommand{\p}{\mathcal{P}}

\newcommand{\cc}{\mathfrak{C}}
\newcommand{\T}{\mathbf{T}}

\begin{abstract}
We show the existence of smoothly embedded closed minimal surfaces in infinite volume hyperbolic $3$-manifolds except some special cases.
\end{abstract}

\maketitle

\section{Introduction}

The existence of minimal surfaces in general $3$-manifolds is one of the classical problems of geometric analysis. There has been numerous foundational results on the existence of closed embedded minimal surfaces in compact $3$-manifolds, or positive curvature case. By \cite{Pi}, every closed $3$-manifold contains a smooth, embedded, closed minimal surface. By geometric measure theory \cite{Fe}, we know that for any compact manifold, there exists an area minimizing surface in every homology class. By \cite{MSY}, it is known that every isotopy class contains a minimal surface in such manifolds. Furthermore, very recently, Song proved Yau's Conjecture and showed that any closed $3$-manifold contains infinitely many minimal surfaces by using the techniques developed by Marques and Neves \cite{So1, IMN}

For noncompact case, the situation is highly different. Even in the constant negative curvature (hyperbolic) case, the answer is unknown.

\begin{question} \label{mainques} Which complete hyperbolic $3$-manifolds contain a smoothly embedded closed minimal surface?	
\end{question} 

In the tree diagram below, green boxes represent known cases, and the red boxes represent the unknown cases for the existence problem.  Note that for some trivial cases like $M=\BH^3$ or $M$ is a Fuchsian Schottky manifold, $M$ contains no closed minimal surface by the maximum principle. In this paper, we give a fairly complete answer to this question by dealing with all red boxes except some special cases listed below.

The summary of the known cases (green boxes) is as follows: For closed hyperbolic $3$-manifolds, Almgren-Pitts min-max theory \cite{Pi} gives a positive answer. For noncompact hyperbolic $3$-manifolds, there are 2 cases: Finite volume, and infinite volume. Recently, several authors proved the existence of minimal surfaces in finite volume noncompact hyperbolic $3$-manifolds \cite{CHMR, HW2, CL, So}, and finished this case. Hence, the only remaining case is the infinite volume hyperbolic $3$-manifolds. 

\

\begin{tikzpicture} [level distance=.7in,sibling distance=12em, 
man/.style={rectangle,draw,fill=green!20},
woman/.style={rectangle,draw,fill=red!20},
edge from parent fork down,
every node/.style = {shape=rectangle, rounded corners=.8ex,	draw, align=center},
grandchild/.style={grow=down,xshift=1em,anchor=west,
	edge from parent path={(\tikzparentnode.south) |- (\tikzchildnode.west)}},
first/.style={level distance=6ex},
second/.style={level distance=12ex}]

\node {Hyperbolic 3-manifolds}
child { node { Infinite Volume  }
	child { node[woman] {Geom. Infinite} 
		child[grandchild,first] {node[woman] {Bounded Geom.}}
		child[grandchild,second] {node[woman] {Unbounded Geom.}}	}
	child { node[man] {Geom. Finite}}
}
child { node[man] { Finite Volume  }
	child[grandchild,first] { node[man] {Closed Manifolds}}
	child[grandchild, second] { node[man] {Cusped Manifolds}}
}

;

\end{tikzpicture}

\

An infinite volume hyperbolic $3$-manifold $M$ contains compact codimension-$0$ submanifold $\C_M$, {\em compact core}, where $M\simeq int(\C_M)$ (See Figure \ref{endsfig}). Then, the ends of $M$ can be thought as the components of $M-\C_M$ (assuming no cusps). Hence, if $\partial \C_M=S_1\cup...\cup S_k$, then each end $\E_i\simeq S_i\times (0,\infty)$ by Marden's Tameness Conjecture \cite{Ag,CG}. The type of such an infinite volume hyperbolic $3$-manifold $M$ is determined by the geometry of these ends. There are two classes: {\em Geometrically finite} and {\em geometrically infinite}. Geometrically infinite ends split into two subclasses depending on the existence of positive injectivity radius bound: \textit{Bounded Geometry} and \textit{Unbounded Geometry}. A basic review on this tree classification of hyperbolic $3$-manifolds is given in Section \ref{sec-HypMan}.

Before stating our main result, we need to define  {\em exceptional manifolds}:

\begin{itemize}
	\item \textit{Type I:} $M$ is a geometrically infinite product manifold.
	\item \textit{Type II:} $M$ is a handlebody (Schottky Manifold).
\end{itemize}

\begin{thm} \label{thm-main}
	Let $M$ be a complete hyperbolic $3$-manifold with finitely generated fundamental group. Then, if $M$ is not an exceptional manifold as listed above, $M$ contains a smoothly embedded closed minimal surface.
\end{thm}

\vspace{.2cm}

\noindent {\em Outline of the Proof:} The main problem in these noncompact manifolds to obtain a closed minimal surface is that the area minimizing sequence of surfaces might escape to infinity, and give an empty limit. The main idea for our construction is to build a barrier near infinity by using the shrinkwrapping surfaces. In particular, shrinkwrapping surfaces give us piecewise-smooth (defective) minimal surfaces. In each end, we construct very special shrinkwrapping surface so that it is either "convex" or "concave" in the end. Then, if we have a convex shrinkwrapping surface in an end, then we can trap area minimizing representative of that end in the compact part of the manifold. If all of the shrinkwrapping surfaces in the ends are concave, then we use min-max methods to obtain a closed embedded minimal surface. For our construction, the regular cusps (rank-2) will not be a problem as our minimal surfaces cannot go very deep in these cusps. Note that throughout the paper, we assume $M$ can only have rank-2 cusps.

\vspace{.2cm}

\noindent {\em Exceptional Manifolds:} As listed above, we have two exceptional families for the existence question. In Type I, $M$ is a product manifold with geometrically infinite ends, i.e. $M\simeq S\times \BR$ where $S$ is a closed surface of genus $\geq 2$.  This case is quite delicate, as one needs to rule out the existence of mean convex foliation in $M\simeq S\times \BR$. For Type II, we have examples like $M=\BH^3$, and Schottky manifolds. For this case, it might be possible to show that there is no closed embedded minimal surfaces in general by constructing CMC foliations. See Sections \ref{sec-Exc} and \ref{sec-remaining} for further discussion.

Organization of the paper is as follows. In Section 2, we give the necessary background, and the relevant results on the problem. In Section 3, we prove our main result. In Section 4, we  discuss the exceptional cases, and further questions, and give some concluding remarks.
 
\subsection{Acknowledgements} We would like to thank Brian Bowditch, Yair Minsky, David Gabai, Antoine Song, Zeno Huang, and Brian White for very valuable conversations. We would like to thank the referees for very valuable comments, and suggestions to improve the paper.

\section{Preliminaries} \label{Prelim}

In this section, we cover the basic definitions and the known results on the problem.

\subsection{Classification of Hyperbolic $3$-Manifolds and Existence Results.} \label {sec-HypMan} \

\vspace{.2cm}

Throughout the paper, we only deal with complete, orientable hyperbolic $3$-manifolds with finitely generared fundamental group. Likewise, all the surfaces will be orientable. To go over the previous results on the existence of minimal surfaces, we follow the tree diagram in the introduction. For more details on the classification of hyperbolic $3$-manifolds, see the survey paper by Minsky \cite{Mi}.

\vspace{.2cm}

\noindent {\bf Finite Volume Hyperbolic $3$-manifolds:} 
These manifold are classified by their compactness:

\vspace{.2cm}

\noindent $\diamond$ {\em Closed Manifolds:} By foundational result of Almgren-Pitts min-max theory \cite{Pi}, every closed Riemannian manifold contains a smoothly embedded minimal surface. The minimal surfaces obtained by this method are \textit{unstable} by construction. See also \cite{MR}. On the other hand, if $M$ is Haken, one can obtain area minimizing surfaces in the homology class of the incompressible surface \cite{MSY,Ha2}. These minimal surfaces are \textit{stable}.
	
	\vspace{.2cm}
	
\noindent $\diamond$ {\em Cusped Manifolds:} In the noncompact case, the finite volume hyperbolic $3$-manifolds is homeomorphic to the interior of a compact $3$-manifold where every boundary component is a torus. The ends of these manifolds are called cusps (rank-2 cusps). These cusps are topologically a solid torus where the core circle is removed. Note that throughout the paper, we will omit rank-1 cusps.
	
Recently, Collin-Hauswirth-Mazet-Rosenberg showed that all these manifolds contain a minimal surface \cite{CHMR} by using min-max methods. They prove that minimal surfaces cannot go very deep in the cusps, hence they obtain smoothly embedded closed minimal surfaces. Again,  these minimal surfaces are also \textit{unstable}. Huang-Wang obtained a similar result for such manifolds with different techniques \cite{HW2}. Their surfaces are least area in their homotopy class, and hence {\em stable}.
	
Note also that recently, Chambers-Liokumovich obtained a very general existence result for any finite volume noncompact $3$-manifold \cite{CL} by generalizing Almgren-Pitts min-max techniques. Song showed that these manifolds contain infinitely many closed embedded minimal surfaces \cite{So}.

\vspace{.2cm}

\noindent {\bf Infinite Volume Hyperbolic $3$-manifolds:} These manifolds are classified by the geometry of their ends \cite{Th}. 

\begin{figure}[h]
	
	\relabelbox  {\epsfxsize=3.5in
		
		\centerline{\epsfbox{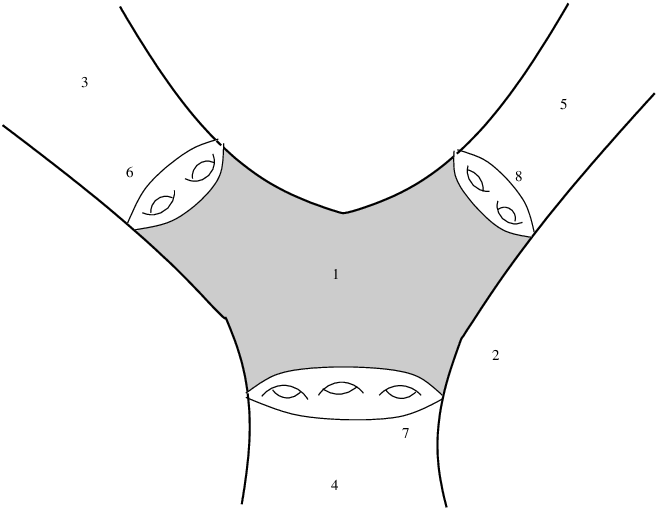}}}

	\relabel{1}{\footnotesize $\C_M$} \relabel{2}{\footnotesize $M$}  \relabel{3}{\footnotesize $\E_1$} \relabel{4}{\footnotesize $\E_2$} \relabel{5}{\footnotesize $\E_3$} 
	\relabel{6}{\scriptsize $S_1$} \relabel{7}{\scriptsize $S_2$} \relabel{8}{\scriptsize $S_3$}

	\endrelabelbox
	
	\caption{\label{endsfig} \small $M$ is an infite volume hyperbolic $3$-manifold with $3$ ends. The shaded region is the compact core $\C_M$.}
	
\end{figure}

Let $M$ be a complete hyperbolic $3$-manifold with infinite volume. Let $\C_M$ be a {\em compact core} of $M$. In particular, $\C_M$ is a compact codimension-$0$ subset of $M$ where $M$ is homeomorphic to interior of $X$. Then assuming no cusps, $M-\C_M=\bigcup_{i=1}^N \E_i$ where $\E_i\simeq S_i\times(0,\infty)$ for some component $S_i$ of $\partial \C_M$ which is a closed surface of genus $\geq 2$ \cite{Ag,CG}. Here, we call $\E_i$ {\em an end of $M$} (See Figure \ref{endsfig}). The geometry of the ends $\{\E_1, \E_2,...,\E_n\}$ describes how complicated the manifold $M$ is. Note that throughout the paper, if $M$ has cusps, we do not treat the cusp regions as an end of $M$, but we call them cusps.

An end $\E_i$ is {\em geometrically finite end} if we can take the boundary component $S_i\subset \partial \C_M$ convex, i.e. $\C_M$ is in the convex side of $S_i$. Geometrically finite ends have very simple geometry, which can be foliated by equidistant surfaces to $S_i$. In particular, if $S^t_i$ is the $t$-equidistant surface to $S_i$ in $\E_i$, then $\{S^t_i\}$ foliates $\E_i$ where $S^t_i$ is a convex surface whose area is increasing exponentially in $t$. 

If an end is not geometrically finite, we call it {\em geometrically infinite end}. In particular, let $\cc(M)$ be the convex core of $M$ which is the smallest convex set in $M$ where the inclusion $\cc(M)\hookrightarrow M$ is a homotopy equivalence. In other words, $\E_i$ is geometrically infinite end in $M$ if $\E_i\subset \cc(M)$. Hence, every component of $M-\cc(M)$ is a geometrically finite end.

Geometrically infinite ends are classified into two categories as follows: Ignoring the cusps, if the geometrically infinite end $\E$ has positive injectivity radius, then we call $\E$ has {\em bounded geometry}. Otherwise, we say that $\E$ has {\em unbounded geometry}. In both cases, there exists a sequence of geodesics $\{\alpha_n\}$ in the end $\E$ escaping to infinity. In other words, if $\E\simeq S\times [0,\infty)$, then for any $K>0$, there exists an $N$, such that $\alpha_n\in S\times [K,\infty)$ for any $n\geq N$ \cite{Bo}.


Now, we are ready to list the results in complete hyperbolic $3$-manifolds with infinite volume.

\vspace{.2cm}

\noindent $\diamond$ {\bf Geometrically Finite Manifolds:} A complete hyperbolic $3$-manifold $M$ is called {\em geometrically finite} if the convex core $\cc(M)$ has finite volume. In other words, every end of $M$ is a geometrically finite end. In this case, assuming nontrivial $H_2(M)$, the existence of a minimal surface in $M$ is straightforward as follows. Consider the area minimizing surface $T$ in the nontrivial homology class of $\xi\in H_2(M)$ \cite{Ha2,Fe}. As any area minimizing surface in $M$ must stay in the convex core $\cc(M)$, $T \subset \cc(M)$. Since $\cc(M)$ is compact, $T$ cannot escape to infinity. Hence, there exists a smoothly embedded stable minimal surface $T$ in $M$. Quasi-Fuchsian Manifolds are examples of geometrically finite manifolds. Note that one can also use the solutions of asymptotic Plateau problem in the universal cover $\BH^3$ in order to get least area surface in such $M$ \cite{An,Co2}. Notice that in this case, we are excluding some Type II Exceptional Manifolds, e.g. Schottky Manifolds (See Section \ref{sec-Exc}). 
	
\vspace{.2cm}

\noindent $\diamond$ {\bf Geometrically Infinite Manifolds with Bounded Geometry:} We call a complete hyperbolic $3$-manifold $M$ is geometrically infinite with bounded geometry if $M$ has a positive injectivity radius bound, and a geometrically infinite end. In other words, ignoring the cusps, there exists $\rho_0>0$ such that the shortest geodesic $\beta$ in $M$ has length $|\beta|>\rho_0$. In particular, $M$ has a geometrically infinite end, and every such end has bounded geometry. The infinite cover of a mapping torus with pseudo-Anosov monodromy (Cannon-Thurston manifolds) are examples of such manifolds \cite{CT}. There is no general existence result for this case.

\vspace{.2cm}

\noindent $\diamond$ {\bf Geometrically Infinite Manifolds with Unbounded Geometry:} We call a complete hyperbolic $3$-manifold $M$ is geometrically infinite with unbounded geometry if $M$ has a geometrically infinite end with unbounded geometry. In particular, $M$ has a sequence of geodesics $\{\alpha_n\}$ which escapes every compact subset of $M$ and $|\alpha_n|\searrow 0$. There is no general existence result for this case, either.

We see that when a complete hyperbolic $3$-manifold $M$ is finite volume, or geometrically finite, $M$ contains a closed embedded minimal surface. Hence, the only remaining cases for Question \ref{mainques} are  infinite volume hyperbolic $3$-manifolds with geometrically infinite ends.

\subsection{Bounded Diameter Lemmas.} \label{bddsec} \

\vspace{.2cm}

In this part, we give three key lemmas about diameters of closed minimal surfaces in hyperbolic $3$-manifolds. These lemmas are crucial for our methods to capture the minimal surfaces in the compact part.

First, we show that diameter of a closed incompressible minimal surface $S$ in a hyperbolic $3$-manifold $M$ can be bounded by a constant depending on its genus $g$, and the injectivity radius of $M$. 

\begin{lem} \label{boundiamlem} [Bounded Diameter Lemma]  
Let $M$ be a hyperbolic $3$-manifold. Let $\rho>0$ be the injectivity radius of $M$.  If $S$ is an incompressible minimal surface in $M$, then there is a constant $C = C(g,\rho)\in \BR^+$ such that  $\mbox{diam}(S) \leq C$. 	
\end{lem}

\begin{pf} As $S$ is incompressible, and the injectivity radius of the manifold is $\rho>0$ for any point $q\in S\cap (M-M_{[0,\rho)})$, $S$ contains an embedded $\rho$-disk $\B_\rho(q)$. Notice that by Gauss' Theorem, as $S$ is a minimal surface, $K_S<K_M=-1$ where $K_S$ is the sectional curvature of $S$, and $K_M$ is the ambient sectional curvature. Then, we have $-\int_S K_S\geq \mbox{area}(S)$. Then, by Gauss-Bonnet, $\int_S K_S=2\pi \chi(S) - \int_{\partial S}\kappa \ dl$ where $\kappa$ is the geodesic curvature. Therefore, as  $S$ is a closed minimal surface in $M$, $\mbox{area}(S)\leq -2\pi\chi(S)= 2\pi(2g-2)$  where $g$ is the genus of the surface. Furthermore, we have $\mbox{area}(B_\rho(q))\geq 2\pi(\cosh(\rho)-1) \geq \pi\rho^2$. 
	
Now, let $\mathcal{P}$ be the maximal collection of points in $S$ which are $2\rho$ apart, and let $N_\rho$ be the number of points in $\mathcal{P}$. Then, we get $N_\rho\cdot \pi \rho^2\leq 2\pi(2g-2)$, which gives an upper bound for $N_\rho$, i.e. $N_\rho\leq \dfrac{4(g-1)}{\rho^2}$. Hence, we have $\mbox{diam}(S)\leq N_\rho\cdot 2\rho=\dfrac{8(g-1)}{\rho}$. The proof follows.
\end{pf}

By \cite[Lemma 1.15]{CG}, this result applies to shrinkwrapped surfaces, too. 

Now, we modify the proof of the above lemma for area minimizing surfaces in a homology class of $M$. Notice that the main difference is that the area minimizing representative may not be incompressible. The main ingredients are Lackenby's general bounded diameter result on closed minimal surfaces in general $3$-manifolds \cite{La}, and genus bound for stable minimal surfaces in hyperbolic $3$-manifolds \cite{BD}.

\begin{lem} \label{boundiamlem2} Let $M$ be a hyperbolic $3$-manifold with $inj(M)=\rho_0>0$. Let $S$ be an embedded closed surface with nontrivial homology, and let $\xi$ be a homology class of $S$ in $M$. If $\Sigma$ is an area minimizing surface in $\xi$, then there exist a constant $C(\xi,\rho_0)>0$ such that   $\mbox{diam}(\Sigma_i) \leq C(\xi,\rho_0)$.
\end{lem}

\begin{pf} By assumption, $\mbox{area}(\Sigma)\leq\mbox{area}(S)$. As $\Sigma$ is area minimizing in $\xi$, it is also stable. Then, by \cite[Lemma 3.3]{BD}, we have $\pi|\chi(\Sigma)|\leq \mbox{area}(\Sigma) \leq 2\pi|\chi(\Sigma)|$. Hence,  $\pi|\chi(\Sigma)|\leq \mbox{area}(\Sigma)\leq\mbox{area}(S)$ gives an upper bound for $|\chi(\Sigma)|$ and the genus of $\Sigma$.
	
	Furthermore, by \cite[Proposition 6.1]{La}, $\mbox{diam}(\Sigma)<|\chi(\Sigma)|f(\kappa,\rho)$ where $\rho$ is the lower bound for injectivity radius, and $\kappa<0$ is the upper bound for the sectional curvature of $M$. By above, we have an upper bound for $\chi(\Sigma)$. Hence, we have $\mbox{diam}(\Sigma) \leq C=|\chi(\Sigma)|f(\kappa,\rho)$. As $\kappa=-1$, and $\rho=\rho_0$ in our case, the constant $C$ only depends on the homology class $\xi$, and $\rho_0$, i.e. $C=C(\xi,\rho_0)$. The proof follows.
\end{pf}

When we have a lower injectivity radius bound, the above lemmas give a uniform diameter bound independent of any region in the manifold. However, when we have a geometrically infinite end with unbounded geometry, because of the short geodesics, we do not have such a such an injectivity radius bound. In this case, we can still have a "pointwise" diameter bound for least area surfaces.

\begin{lem} \label{boundiamlem3} [Pointwise Diameter Bound for Ends with Unbounded Geometry]
Let $M$ be a hyperbolic $3$-manifold, and let $\E\simeq T\times[0,\infty)$ be a geometrically infinite end with unbounded geometry. let $K$ be a compact set in $\E$. Let $\Y_K$ be space of closed surfaces in $\E$ in the homotopy type of $T$ and with nontrivial intersection with $K$. Let $\Sigma$ be the least area surface in $\Y_K$. Then, there exists $C_K>0$ such that $d_K(\Sigma)=\sup\{d(x,K) \mid x\in \Sigma\} <C_K$.	
\end{lem}

\begin{pf} Recall that $\Y_K=\{S\subset \E \mid S\sim T \mbox{ and } S\cap K\neq \emptyset\}$. Let $\Sigma$ be the least area surface in $\E$. Because of the nontrivial intersection condition, $\Sigma$ will be piecewise smooth, where the nonsmooth parts are in $\Sigma\cap K$. 
	
Let $\E_{thick}$ be the thick part of $\E$. In particular, $\E_{thick}=\E\cap M_{thick}$ is connected, and its injectivity radius is greater than the Margulis constant $\e_0$. Let $\E_{thin}=E-E_{thick}$ be the thin part of $\E$, which is the union of Margulis tubes of short geodesics in $\E$. 

Let $\Sigma_{thick}=\Sigma\cap \E_{thick}$ and $\Sigma_{thin}=\Sigma\cap \E_{thin}$. Notice that $\Sigma_{thin}$ is union of very long and thin subsurfaces in Margulis tubes where we cannot bound the diameter as Margulis tubes can have very large diameters. However, we can bound the number of components in $\Sigma_{thin}$ as follows.  Let $m$ be the number of components of $\Sigma_{thick}$. Let $\{V_i\}$ be the set of components of $\Sigma_{thick}$, i.e. $\Sigma_{thick}=\bigcup_{i=1}^m V_i$. Then, because of the injectivity radius bound $\e_0>0$ on $\Sigma_{thick}$, the area of the each component must be at least the area of one disk, i.e. $|V_i|>2\pi(\cosh(\e_0)-1)>\pi\e_0^2$ by Lemma \ref{boundiamlem}. By Gauss-Bonnet, total area of $\Sigma$ is uniformly bounded, i.e. $|\Sigma|<2\pi(2g-2)$ (Lemma \ref{boundiamlem}). This implies $m\cdot(\pi\e_0^2)<2\pi(2g-2)$. Hence, we have $m<\dfrac{2\pi(2g-2)}{\pi\e_0^2}$.

This implies the number of components $m$ in $\Sigma_{thick}$ is only depends on the genus of the surface $S$, i.e. ($m$ depends only on $\E$). Now, we assumed that $\Sigma\cap K\neq \emptyset$. Now, we claim that $d_K(\Sigma)=\sup\{d(x,K) \mid x\in \Sigma\}$ can be bounded by the sum of the diameters of Margulis tubes near $K$. 

Let’s start with an easy case. Assume $K$ is in the thick part for simplicity. If m=2, then we get the bound as follows: Let $C_0$ be diameter bound for the thick part. Let $D_1$ be the diameter of largest Margulis tube whose distance to $K$ is less than $C_0$. Then distance $d_K(S)$ of $S$ to $K$ could be at most $C_0 + D_1$ as there can be at most 2 components in the thick part. Now, we iterate this argument for general $m$.

Let $T_1,.., T_m$ be the components of $\Sigma_{thick}$. Recall that $\Sigma$ is a connected surface. Hence, for each $i$, there is a $j$ such that $T_i$ and $T_j$ are connected by a thin surface in a Margulis tube, and $d(T_i,T_j)$ is bounded by the diameter of the tube. Now, order the components $\{T_i\}$ by their distances to the compact set $K$. Without loss of generality, assume $T_1\cap K\neq 0$. So, $d(K,T_1)=0$.

Now, let $D_1$ be the diameter of largest Margulis tube whose distance to $K$ is less than $C_0$.  This implies $d(T_2,K)<C_0+D_1$.
Let $D_2$ be the diameter of largest Margulis tube whose distance to $K$ is less than  $C_0+D_1$. This implies $d(T_3,K)<C_0+D_1+ D_2$.
Let $D_{i+1}$ be the diameter of largest Margulis tube whose distance to $K$ is less than $C_0+D_1+D_2+…+D_i$. $d(T_{i+1},K)<C_0+D_1+...+D_i$

This gives $d(T_m,K)<C_0+D_1+...+D_{m-1}$. As $T_m$ is the farthest component, define $C_K= 2C_0+ D_1+...+D_m$. As the diameters of the components are bounded, $d_K(\Sigma) <C_K$. The proof follows.
\end{pf}

\begin{rmk} Intuitively, our least area surface can make huge jumps through the "Margulis tubes" in the end. In the lemma above, we bound the number of jumps the surface can make, and hence by fixing its starting point, we give a bound for its diameter. 
	
Notice that first two bounded diameter lemmas above are only valid when there exist a positive injectivity radius bound, and it can be applied anywhere in the manifold independent of the location in $M$. However, the last one gives a "pointwise" diameter bound, where one must specify a compact set $K\subset M$ first to use this lemma. The advantage of this lemma is that it does not need a positive injectivity radius bound.
\end{rmk}

\subsection{Cusps and Margulis Tubes.} \

\vspace{.2cm}

In this part, we go over the basic definitions, and relevant results on cusps and short geodesics in hyperbolic $3$-manifolds. For further details, see \cite{BP, Ma}. Note that throughout the paper, we omit rank-$1$ cusps, and assume all the cusps are rank-$2$. See Section \ref{sec-remaining} for discussion on rank-$1$ cusps.

\begin{figure}[b]
	
	\relabelbox  {\epsfxsize=3.5in
		
		\centerline{\epsfbox{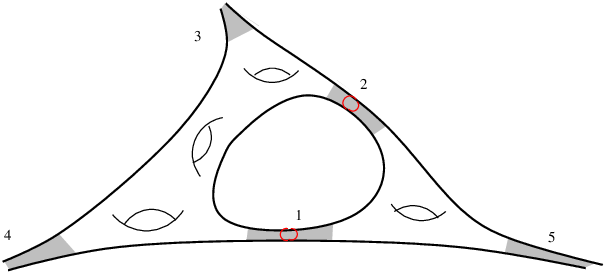}}}

	\relabel{1}{\footnotesize $\gamma_1$} \relabel{2}{\footnotesize $\gamma_2$}  \relabel{3}{\footnotesize $X_1$} \relabel{4}{\footnotesize $X_2$} \relabel{5}{\footnotesize $X_3$}

	\endrelabelbox
	
	\caption{\label{cuspfig} \small $2$-dimensional analogous picture for Margulis Tubes and Cusps. $\gamma_1$ and $\gamma_2$ short geodesics, and shaded regions represent Margulis Tubes and Cusps.}
	
\end{figure}

Any finite volume, complete, non-compact hyperbolic $3$-manifold $M$ is homeomorphic to the interior of a compact $3$-manifold $N$ where $\partial N$ is a collection of tori. In particular, the ends of a finite volume, complete hyperbolic $3$-manifolds are homeomorphic to $T^2\times[0,\infty)$. One can think of these ends as a solid torus where the core circle is removed. Such ends are called {\em cusps}. Informally, one can use the analogy between a punctured surface with hyperbolic structure, where the end has $S^1\times[0,\infty)$ structure near the puncture (See Figure \ref{cuspfig}). 

The geometry of the cusps are well understood, as they are generated by parabolic isometries in the representation of the fundamental group of the manifold. In particular, every cusp is foliated by constant mean curvature 1 tori, induced by the corresponding horospheres in the universal cover $\BH^3$.

On the other hand, a similar structure is also true for neighborhoods of sufficiently short geodesics in $M$. Let $\epsilon_3$ be the Margulis constant for hyperbolic $3$-manifolds, which is known to be between $0.1$ and $0.6$.
Let $M_{[0,\epsilon_3)}$ is the part of $M$ where the injectivity radius is smaller than $\e_3$. We call $M_{[0,\epsilon_3)}$ the {\em thin part} of $M$. By Margulis Lemma, the structure of the thin parts of the hyperbolic $3$-manifold is well understood. It is either homeomorphic to $T^2\times[0,\infty)$ (cusp), or homeomorphic to a solid torus where the core circle is a short geodesic (See Figure \ref{cuspfig} for an analogous picture). These solid torus neighborhoods $\N_r(\gamma)$ of short geodesics $\gamma$ are called {\em Margulis Tubes}. The diameter of the maximal tubes depends on the length of the geodesic. In particular, if $l=\mbox{length}(\gamma)\searrow 0$, then $r\nearrow \infty$ (roughly $r\sim -\log l$). Their geometric structure is very similar to the cusps. In a way, one can think of cusps as Margulis tubes for "length $0$ geodesics".

Note that cusps and short geodesics can appear in not only finite volume hyperbolic $3$-manifolds, but also infinite volume hyperbolic $3$-manifolds. Even though cusps can topologically be considered as an end of the $3$-manifold, because of their very special structure, we do not call them an end of the manifold, instead we call them {\em cusps}. Hence, when we say $\E$ is an end  of a hyperbolic $3$-manifold $M$, we always mean $\E\simeq S\times [0,\infty)$ where $S$ is a closed surface of genus $\geq 2$.

Note also that even for a complete hyperbolic 3-manifold $M$ with finitely generated fundamental group, there can be infinitely many components in the thin part of $M$ for every positive $\e$ smaller than the Margulis constant \cite{BO}. However, Sullivan proved that a complete hyperbolic 3-manifold $M$ with finitely generated fundamental group has only finitely many cusps \cite{Su}. So, by choosing the compact core $\C_M$ accordingly, we can assume that there is no cusp region in the ends, i.e. $\E_i\simeq S_i\times [0,\infty)$

In \cite{CHMR}, the authors showed the following strong result about minimal surfaces in cusp neighborhoods when proving the existence of minimal surfaces in finite volume hyperbolic $3$-manifolds. The following lemma states that minimal surfaces in a hyperbolic $3$-manifold cannot go very deep in a cusp.

\begin{lem} \label{cusplem} \cite[Proposition 8]{CHMR} \cite{CHMR2} 
Let $M$ be a complete hyperbolic $3$-manifold, and let $X=T^2\times[0,\infty)$ be a cusp in $M$. Let $\Sigma$ be a closed minimal surface in $M$ whose Morse index is bounded below by $k$. Then, there exists $C_X(k)>0$ such	that $\Sigma\cap X\subset T^2\times[0, C_X(k)]$.	
\end{lem}

Note that all the surfaces we construct in this paper will be stable except where we use Corollary \ref{minmaxcor} where we have a index bound. Therefore, the lemma above will suffice for our purposes. For different versions of the lemma above, see also \cite[Theorem 5.9]{HW2} and \cite[Lemma 2.5]{Ha}.

\begin{rmk} \label{rmk-tubelem} [Tube Lemma] In an earlier version of this paper, we used Tube Lemma in our construction. The statement was basically a least area surface cannot come very close to the core circle of a Margulis tube \cite{HW,Ha}. Thanks to the referee and Brian Bowditch, we realized a technical problem in our generalization, and we removed the  tube lemma from our construction. In particular, we had used this lemma in order to show the existence of closed embedded minimal surfaces in the geometrically infinite ends with unbounded geometry. Instead of using tube lemma, by using Lemma \ref{boundiamlem3}, we adapt our shrinkwrapping technique to the ends with unbounded geometry to resolve this case.
\end{rmk}

On the other hand, the tube lemma indeed may not be true in general. The main problem is that a minimal surface may intersect a Margulis Tube $\N_r(\gamma)$ in an essential, thin long annulus $\A$ where $\partial \A=\alpha^+\cup\alpha^-$ are homotopic to the longitude of the Margulis tube. In such a case, if the length of the geodesics $l_n=|\gamma_n|\searrow 0$, then the area of such thin annuli $|\A_n|=l_n.\sinh{r_n}<\sqrt{\frac{\sqrt{3}l_n}{4\pi}}\to 0$  \cite[P 1315]{HW}. So, the area comparison do not work in such situation. However, if we have a complex condition like $\frac{|\theta_n|}{l_n}>\sqrt[4]{3\pi^2}$, then Huang-Wang gets the necessary estimates to prove the tube lemma, i.e. a least area surface stays away from the core of the Margulis tube. For further details, see \cite[Section 3.1]{HW}. We would like to thank the referee for pointing out this problem.

Following heuristic construction by Brian Bowditch shows how subtle the problem can be.

\vspace{.2cm}

\noindent {\em Bowditch's Tube Example:} Let $\alpha,\beta$ be two disjoint simple closed curves on a closed surface $S$, which cut $S$ into three subsurfaces, $B,A,C$ with $A$ in the middle. Let $f_1, f_2, ..., f_6$ be pseudo-Anosov maps defined respectively on the surfaces $B, A, C, A \cup B, A$, and $A \cup C$. Let $F_n$ be the composition $f_1^n, f_2^n, ... ,f_6^n$. This will be pseudo-Anosov map on $S$ for all sufficiently large $n$. Let $M_n$ be the mapping torus of $F_n$.

\begin{wrapfigure}{r}{1.5in}
	\vspace{-.2cm}
	\relabelbox  {\epsfxsize=1.5in
		
		\centerline{\epsfbox{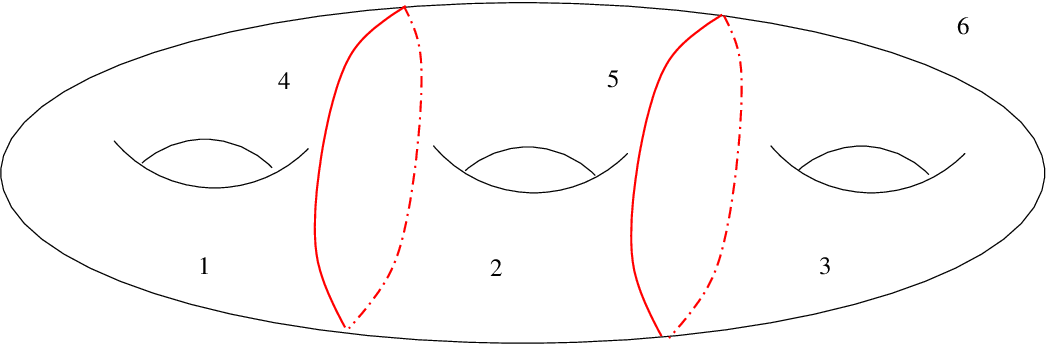}}}
	\relabel{1}{\tiny $B$}
	\relabel{2}{\tiny $A$}
	\relabel{3}{\tiny $C$}
	\relabel{4}{\tiny $\alpha$}
	\relabel{5}{\tiny $\beta$}
	\relabel{6}{\scriptsize $S$}
	\endrelabelbox
	\vspace{-.2cm}
\end{wrapfigure}

Given any $\e>0$, for sufficiently large $n$, $M_n$ will decompose into two $\e$-Margulis tubes (about $\alpha$ and $\beta$), together with {\em bands}, topologically of the form $F \times I$, where $F$ is one of the base surfaces, $A,B,C, A \cup B$ or $A \cup C$, and $I$ is a real interval.  These meet along their boundaries, and one can draw explicitly the combinatorial picture.

By choosing $n$ large enough, these bands can each be made arbitrarily long in the $I$-direction.  This means that, for any given $d$, we can ensure that any fiber of $M_n$ disjoint from the Margulis tubes must have intrinsic diameter at least $d$, since it has to cross at least one of the bands (In fact, one only needs to consider the two bands with base surface arising from $A$).  If $d$ is large enough in relation to $\e$, then this fiber cannot be a minimal surface, since such a surface has intrinsic curvature at most $-1$.

On the other hand, there must be a least area surface homotopic in $M_n$ to the fiber by \cite{MSY}. Lifting the infinite cyclic cover of $M_n$ gives a least area surface in a doubly degenerate manifold $\wt{M}_n$. Such a least area surface can come very close to the core circles of Margulis tubes of $\alpha$ and $\beta$.

\subsection{Shrinkwrapping.} \label{shrinksec} \

\vspace{.2cm}

Throughout the paper, we will frequently use shrinkwrapping technique to obtain "defective" minimal surfaces. The shrinkwrapping technique was introduced by Calegari and Gabai in \cite{CG} where they prove the celebrated {\em Marden's Tameness Conjecture}. A simpler version of this technique can also be seen at \cite{Co1}.

\begin{figure}[h]
	
	\relabelbox  {\epsfxsize=3.5in
		
		\centerline{\epsfbox{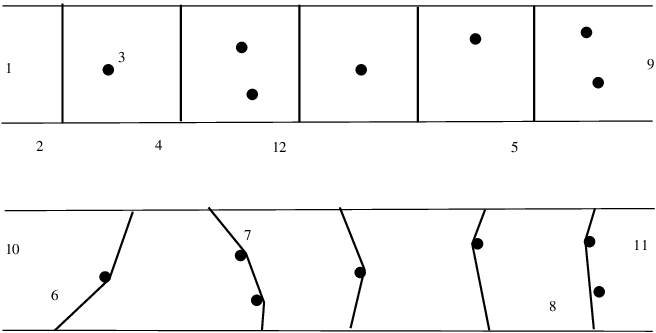}}}

	\relabel{1}{\footnotesize $M$} \relabel{2}{\footnotesize $S\times 0$}  \relabel{3}{\footnotesize $\alpha_1$} \relabel{4}{\footnotesize $S\times 1$} \relabel{5}{\footnotesize $S\times m$} 
	\relabel{6}{\footnotesize $\Sigma_1$} \relabel{7}{\footnotesize $\Sigma_2$}  \relabel{8}{\footnotesize $\Sigma_m$} \relabel{9}{\footnotesize $\E$} 	
	\relabel{10}{\footnotesize $M$} \relabel{11}{\footnotesize $\E$} 
	\relabel{12}{\footnotesize $S\times 2$} 	
	\endrelabelbox
	
	\caption{\label{shrinkfig} \small By shrinkwrapping, we deform $S\times k$ into a $\Gamma$-minimal surface $\Sigma_k$ in the end $\E$.}
	
\end{figure}

In particular, shrinkwrapping is "tightening" a given "incompressible surface" $S$ outside of a finite collection of geodesic curves $\Gamma$. At the end of this "tightening" process, one obtains an isotopy of $S$ into a "defective minimal surface" $\Sigma$, where $\Sigma$ is minimal except along some curve segments along $\Sigma\cap \Gamma$ (See Figure \ref{shrinkfig}). We call such a surface $\Gamma$-minimal, i.e. $\Sigma$ is smooth with mean curvature $0$ everywhere except $\Sigma\cap\Gamma$. The geometry of the surface near these singular segments is also well-understood along the lines of {\em thin obstacle problem}. Here is the main technical result of \cite{CG}. Note that $S$ is $2$-incompressible relative to $\Gamma$ means that any essential compressing disk for $S$ intersects $\Gamma$ at least twice. 

\begin{lem}  \label{shrinklem} [Shrinkwrapping] \cite[Theorem 1.10]{CG}
Let $M$ be a complete, orientable hyperbolic $3$-manifold, and let $\Gamma$ be a finite collection of pairwise disjoint simple closed geodesics in $M$. Furthermore, let $S \subset M-\Gamma$ be a closed embedded $2$-incompressible surface rel. $\Gamma$ which separates some component of $\Gamma$ from another. Then S is homotopic to a $\Gamma$-minimal surface $T$ via a homotopy $F : S \times[0, 1] \to M$ such that

\begin{enumerate}
\item $F(S \times 0) = S$,
\item $F(S \times t) = S_t$ is an embedding disjoint from $\Gamma$ for $0\leq t< 1$,	
\item $F(S \times 1) = T$,
\item if $T'$ is any other surface with these properties, then $\mbox{area}(T) \leq \mbox{area}(T')$
\end{enumerate}
	
\end{lem}

Note that Calegari-Gabai used this lemma to prove celebrated Marden's conjecture: Any hyperbolic $3$-manifold is homeomorphic to the interior of a compact $3$-manifold with boundary. In other words, any end $\E$ of a hyperbolic $3$-manifold has the product structure $\E\simeq S\times [0,\infty)$ for some closed surface $S$ of genus $\geq 2$. 

Calegari-Gabai needed this complicated $2$-incompressibility condition on $S$ and $\Gamma$ as a priori they did not know the product structure of the end. However, with the help of Marden's Conjecture, we can simplify the statement of shrinkwrapping lemma by using the product structure of the end. In particular, the set $\Gamma$ can be chosen as only two geodesics which are sufficiently deep in the end, and in different sides of a level set of the end as follows:

\begin{cor} \label{cor-shrink} Let $M$ be a hyperbolic $3$-manifold, and $\E\simeq S\times [0,\infty)$ be a geometrically infinite end of $M$. Let $\alpha_1$ and $\alpha_2$ be two geodesics in $\E$ where
	
\begin{itemize}
	\item $\alpha_1$ and $\alpha_2$ are separated by $S_0=S\times\{t_0\}$ for some $t_0>0$.
	\item $\alpha_1$ and $\alpha_2$ are sufficiently deep in the end $\E$, i.e. $d(\alpha_i,S\times\{0\})>C_0$ where $C_0$ is the bounded diameter constant given either by Lemma \ref{boundiamlem} (when $\E$ has bounded geometry), or by Lemma \ref{boundiamlem3} (when $\E$ has unbounded geometry).
\end{itemize} 

Then, $S_0$ is homotopic to a piecewise smooth minimal surface $\Sigma$ where the only nonsmooth parts are in $\Sigma\cap (\alpha_1\cup\alpha_2)$.	
\end{cor}

\begin{pf} By using Lemma \ref{shrinklem}, all we need to show that $\Gamma=\{\alpha_1,\alpha_2\}$ satisfy the shrinkwrapping conditions. Notice that $S_0=S\times \{t_0\}$ is already incompressible in $\E$. Therefore there is no compressing essential disk in $E$. Furthermore, $S_0$ separates $\alpha_1$ and $\alpha_2$. So, if we can show that the $\Gamma$-shrinkwrapping surface $\Sigma$ completely stays in $\E$, we are done. In order to make sure that, we will choose our geodesics sufficiently deep in the end by using bounded diameter lemmas.
	
Recall that the $\Gamma$-shrinkwrapping surface $\Sigma$ is the least area surface among the surfaces separating $\alpha_1$ and $\alpha_2$, and homotopic to $S_0$. If $\E$ has bounded geometry, let $C_0$ be the constant given by Lemma \ref{boundiamlem}. If $\E$ has unbounded geometry, let $C_0$ be the constant for the compact set $K=S\times\{0\}$ given by Lemma \ref{boundiamlem3}. As $d(\alpha_i,S\times\{0\})>C_0$, and $\Sigma$ must be between $\alpha_1$ and $\alpha_2$, we have $\Sigma\subset S\times (0,\infty)\in \E$. Hence, $\Gamma=\{\alpha_1,\alpha_2\}$ and $S_0$ satisfy shrinkwrapping conditions in $\E$. Then, by Lemma \ref{shrinklem}, $\Sigma$ is a $\Gamma$-minimal surface, i.e. smooth minimal surface everywhere except where it touches its barriers $\alpha_1$ and $\alpha_2$. The proof follows.
\end{pf}

Before finishing our discussion on shrinkwrapping surfaces, we need to give the regularity results near their "defective" parts. As the lemma above describes, the surface $\Sigma$ is minimal except the {\em coincidence set} $L:=\Sigma\cap\Gamma$. The behavior of $\Sigma$ near $L$ is well-understood, and it is called {\em Thin Obstacle Problem} in the literature. The following lemma summarizes these regularity results. For further details, see \cite[Section 1.7]{CG}.

\begin{lem} \label{thinlem} [Thin Obstacle Problem] \cite[Lemma 1.31]{CG}
	Let $\Sigma$ be a $\Gamma$-minimal surface in $M$ as described above. Then, there exists a parametrization of $\Sigma$ with $u:S\to M$   such that the derivative $du$ along local sheets of $\Sigma$ is continuous from each side along the coincidence set $L$, and continuous at noninterior points.
\end{lem}

\begin{rmk} Note that Calegari-Gabai needed a finer understanding of $\Sigma$ near coincidence set to prove $\Sigma$ is indeed a $CAT(-1)$ surface. In this paper, we only need that $\Sigma$ has tangent planes in both sides of coincidence sets, which follows from the lemma above.	\end{rmk}

\subsection{Min-Max Surfaces in Non-Compact Manifolds.} \

\vspace{.2cm}

In this section, we give a recent min-max result which is the key ingredient to finish off a very important case in our main result. By using Almgren-Pitts Min-Max Theory, Montezuma and Song independently constructed embedded minimal surfaces in non-compact manifolds. Their result is more general, but we give the adapted version here.

\begin{lem} \label{minmaxlem} \cite{Mo, So}
	Let $M$ be a complete non-compact $3$-manifold. If $M$ contains a bounded open set $\Omega$ such that $\overline{\Omega}$ is a smooth strictly mean-concave boundary, then there exists a finite area embedded minimal surface $\Sigma$ in $M$ intersecting $\Omega$.		
\end{lem}

\begin{rmk} \label{minmaxrmk} Note that Song's recent result \cite[Corollary 5]{So} is more general, and indeed we are using his version here.  Montezuma proved a very similar result in \cite{Mo}, but it does not apply to the case when $M$ has a cusp. 
\end{rmk}

Notice that if there is no cusp in $M$, then the finite area minimal surface is automatically closed. However, if $M$ has a cusp, the finite area minimal surface may not be closed, as it can have a non-compact part in the cusp. However, hyperbolic $3$-manifolds and cusps are special. By combining Lemma \ref{cusplem} and above lemma, we have the following stronger result for hyperbolic $3$-manifolds.

\begin{cor} \label{minmaxcor} 
Let $M$ be a complete hyperbolic $3$-manifold. If $M$ contains a bounded open set $\Omega$ such that $\overline{\Omega}$ is a smooth strictly mean-concave boundary, then there exists a closed embedded minimal surface $\Sigma$ in $M$ intersecting $\Omega$.		
\end{cor}

\begin{pf} To prove this corollary, we will look into Song's result above in the hyperbolic $3$-manifold case. In the proof of Theorem 4 and Corollary 5 in \cite{So}, Song obtains the minimal surface as a limit of closed minimal surfaces in larger and larger domains exhausting the manifold. In particular, let $\{U_i\}$ be sequence of domains exhausting $M$ such that $U_i\subset U_{i+1}$ and $M=\bigcup U_i$. Then, Song obtains a closed minimal surface $S_i$ in $U_i$ in two different cases in \cite[Theorem 4]{So}. Then, the minimal surface mentioned in the Theorem 4 and Corollary 5 in \cite{So} is the limit minimal surface $S=\lim S_i$. Furthermore, the closed minimal surfaces $\{S_i\}$ are stable in one case, and Morse index at most 1 in the other case.
	
By Lemma \ref{cusplem}, when $M$ is hyperbolic, minimal surfaces with bounded index cannot go deep in the cusp regions. That means the closed minimal surfaces in the sequence $\{S_i\}$ cannot go deep in the cusps. Then, the limit finite area minimal surface $S=\lim S_i$ must be a closed surface, too.	
\end{pf}


\subsection{Area Minimizing Surfaces in $3$-Manifolds.} \

\vspace{.2cm}

In this section, we quote famous existence results for area minimizing representatives in a homotopy or homology class. These results will play crucial role in our construction. The first one is the celebrated result of geometric measure theory.

\begin{lem} \label{AMSlem} \cite{Fe,Ha2}  Let $M$ be a closed or mean convex compact orientable $3$-manifold. Let $\xi$ be a nontrivial homology class in $H_2(M)$. Then there exists a smoothly embedded surface $\Sigma$ which has the smallest area among the surfaces in $\xi$.	
\end{lem}

The second result is the existence of least area representative in a homotopy class by Meeks-Simon-Yau.

\begin{lem} \label{LAlem} \cite{MSY} Let $S$ be an embedded surface in a Riemannian $3$-manifold $M$. After series of compressions, isotopies, and collapsings of boundaries of $I$-bundles to their cores, $S$ can be realized as a union (possibly empty) of disjoint embedded minimal surfaces.
\end{lem}

The third one is on least area representatives of incompressible surfaces.

\begin{lem} \label{incomplem} \cite{HS}  Let $M$ be a closed or mean convex compact orientable, irreducible $3$-manifold. Let $S$ be an incompressible surface in $M$. Then, there is a smoothly embedded least area surface $\Sigma$ in the isotopy class of $S$.
\end{lem}

Note that we quoted these results in a way which we are going to use in the remainder of the paper. Some results that were stated for closed manifolds can be applied to compact mean convex manifolds by using maximum principle \cite[Section 6]{HS}.

\

\section{The Main Result}

In this section, we prove our main theorem. Before giving our main result, we need to define {\em exceptional manifolds.} Note that we assume $M$ has no rank-1 cusps. A handlebody is a solid (filled) genus-$g$ surface in the homotopy type of bouquet of $g$ circles.

\vspace{.2cm}

\noindent {\bf Exceptional Manifolds:} 

\begin{itemize}
	\item \textbf{Type I:} $M$ is a geometrically infinite product manifold, i.e. $M\simeq S\times \BR$ for some closed surface of genus $\geq 2$.
	
	\vspace{.2cm}
	
	\item \textbf{Type II:} $M$ is a handlebody (Schottky Manifold).
\end{itemize}

Now, we can state our main result.

\begin{thm} \label{mainthm}
	Let $M$ be a complete hyperbolic $3$-manifold with finitely generated fundamental group. Further assume that $M$ is not an exceptional manifold as listed above. Then, there exists a smoothly embedded, closed minimal surface $\Sigma$ in $M$.
\end{thm}

\noindent {\em Outline of the Proof:} We split the proof into 5 steps. In the first two steps, we construct special shrinkwrapping surfaces in each end. In Step 3, we show that these special shrinkwrapping surfaces will be either convex or concave in the end. Then, in Step 4, we study the case when one of the ends contain a convex shrinkwrapping surface. In Step 5, we use min-max method to show the existence of a closed embedded minimal surface if all shrinkwrapping surfaces in the ends are concave.

\begin{pf} Let $M-\C_M=\bigcup_{i=1}^N \E_i$ where $\C_M$ is the compact core of $M$, and $\{\E_i\}$ are the ends of $M$. As there can only be finitely many cusps, we assume $\E_i\simeq S_i\times [0,\infty)$ where $S_i$ is a closed surface of genus $\geq 2$ in $\partial\C_M$. 
	
In the following two steps, we construct a special shrinkwrapping surface $\Sigma_i$ in each end $\E_i$. These special shrinkwrapping surfaces are obtained by minimizing area among the surfaces separating two sufficiently distant short geodesics in the end (Corollary \ref{cor-shrink}). Being special means the shrinkwrapping surface can touch at most one of these geodesics. If an end is geometrically finite, we take the convex surface in the boundary of the convex core corresponding to that end $\E_i$ as our special surface $\Sigma_i$.

Let $C_1$ be the constant in the first bounded diameter lemma for $M'$ (Lemma \ref{boundiamlem}). Let $C_2$ be the constant in the second bounded diameter lemma for $M'$ (Lemma \ref{boundiamlem2}). Let $C_0=2\max\{C_1,C_2\}$.

\vspace{.2cm}
	
\noindent {\bf Step 1:} Let $\E$ be geometrically infinite end. Then, $\E$ contains a special shrinkwrapping surface $\Sigma$. 

\vspace{.2cm}

\noindent {\em Proof of Step 1:} First, assume that $\E\simeq S\times [0,\infty)$ is a geometrically infinite end with bounded geometry in $M$. Then, there exists $\{\beta_n\}$ be an exiting sequence of geodesics in $\E$. Let $\alpha_1$ be a geodesic in  $\E$ where $d(\alpha_1,S\times\{0\})>C_0$. Assume $\alpha_1\in S\times [0,d_0)$ for some $d_0>0$. Let $\alpha_2$ be another geodesic with $S\times (d_0,\infty)$ such that $d(\alpha_1,\alpha_2)>C_0$. Then, $S_0=S\times\{d_0\}$ separates $\alpha_1$ and $\alpha_2$.

Let $\Sigma$ be the shrinkwrapping surface of $S_0$ with respect to $\Gamma=\{\alpha_1,\alpha_2\}$ in $\E$ given by Corollary \ref{cor-shrink}. Furthermore, as we required $d(\alpha_1,\alpha_2)>C_0$, $\Sigma$ cannot intersect both  $\alpha_1$ and $\alpha_2$. Hence, in the bounded geometry case, we obtain a special shrinkwrapping surface $\Sigma$ in $\E$.

Now, assume that $\E$ is a geometrically infinite end with unbounded geometry in $M$. Since we have arbitrarily short geodesics in $\E$, we cannot use bounded diameter lemmas as in the previous step to obtain our shrinkwrapping surfaces. Instead, we use "pointwise bounded diameter lemma" (Lemma \ref{boundiamlem3}), in order to make sure our shrinkwrapping surface intersect at most one geodesic.

Parametrize $\E\simeq S\times [0,\infty)$. Again, change the metric of $\E$ in $S\times[0,\e)$ so that it becomes mean convex towards inside the end. Let's call this space with modified metric $\E'$.

Let $K=S\times\{0\}$. Let $C_K>0$ be the constant given by Lemma \ref{boundiamlem3}. Let $\alpha_1$ be a geodesic in $\E$ such that $d(\alpha_1,K)>C_K$. 


Let $b=\sup\{t\mid \alpha_1\cap S\times\{t\}\neq \emptyset\}$. Let $\E''=S\times[b+2,\infty)$. $K'=S\times\{b+2\}$. By apply Lemma \ref{boundiamlem3} in $\E''$ for the compact set $K'$, we get another another constant $C_{K'}>0$. Now let $\alpha_2 \in \E''$ such that $d(K',\alpha_2)>C_{K'}$. 

Now, we have $S\times \{b+1\}$ separates $\alpha_1$ and $\alpha_2$. Let $\Sigma$ be the shrinkwrapping surface of $S\times\{b+1\}$ with respect to $\Gamma=\{\alpha_1,\alpha_2\}$ given by Corollary \ref{cor-shrink}. In particular, $\Sigma$ is the least area surface among the surfaces separating $\alpha_1$ and $\alpha_2$. Hence, by the choice of $C_{K'}$, $\Sigma$ cannot intersect both geodesics $\alpha_1$ and $\alpha_2$. Step 1 follows. \hfill $\Box$

\vspace{.2cm}

Now, for each special shrinkwrapping surface $\Sigma$ in $\E$, we have the following three mutually exclusive cases. Let $\alpha_1$ be the curve closer to the compact core, and $\alpha_2$ is the farther one. 
	
\begin{enumerate}
	\item If $\Sigma\cap \alpha_1\neq \emptyset$, we call $\Sigma$ is a {\em convex shrinkwrapping surface}, 
	\item  If $\Sigma\cap \alpha_2\neq \emptyset$, we call $\Sigma$ is a {\em concave shrinkwrapping surface},
	\item  $\Sigma\cap \alpha_i=\emptyset$ for $i=1,2$.
\end{enumerate}

Notice that in case (3), $\Sigma$ is already smoothly embedded closed minimal surface in $M$, and we are done. So, we will omit this case. 

Following step clarifies why we call $\Sigma$ a convex or concave shrinkwrapping surface (See Figure \ref{shrinkfig2}).


\vspace{.2cm}
	
\noindent {\bf Step 2:} If $\Sigma$ is convex shrinkwrapping surface, then the bounded component of $\E-\Sigma$ is mean convex. If $\Sigma$ is concave shrinkwrapping surface, then the unbounded component of $\E-\Sigma$ is mean convex.

\vspace{.2cm}

\noindent {\em Proof of Step 2:} Assume that $\Sigma$ is a concave shrinkwrapping surface, i.e. $\Sigma\cap \alpha_2\neq \emptyset$ . Recall that $\Sigma$ is obtained by shrinkwrapping the surface $S\times\{t_0\}$ in the end $\E$ which separates the geodesics $\alpha_1$ and $\alpha_2$. By the proof of Lemma \ref{shrinklem}, $\Sigma$ is a limit of smooth embedded surfaces $F_n$ which are minimal except a small ($r_n\searrow 0$) neighborhood of geodesic  $\alpha_2$. Let $\tau$ be a proper curve segment connecting $\alpha_2$ and $\alpha_1$. Then, the algebraic intersection number $\langle S_{m_0}, \tau\rangle=1$. By construction, $F_n$ and $\Sigma$ are homotopic to $S_{m_0}$, then we have $\langle \Sigma, \tau\rangle=1$ and $\langle F_n, \tau\rangle=1$ for any $n$ by Lemma \ref{shrinklem}. This implies $\Sigma$ and $F_n$ individually separate $\alpha_1$ and $\alpha_2$ in $\E$.

\begin{figure}[t]
	
	\relabelbox  {\epsfxsize=3.5in
		
		\centerline{\epsfbox{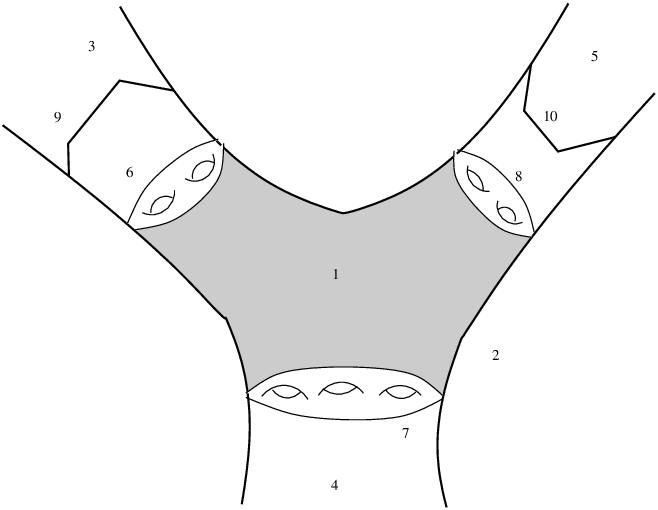}}}

	\relabel{1}{\footnotesize $\C_M$} \relabel{2}{\footnotesize $M$}  \relabel{3}{\footnotesize $\E_1$} \relabel{4}{\footnotesize $\E_2$} \relabel{5}{\footnotesize $\E_3$} 
	\relabel{6}{\scriptsize $S_1$} \relabel{7}{\scriptsize $S_2$} \relabel{8}{\scriptsize $S_3$}
	\relabel{9}{\scriptsize $\Sigma_1$} \relabel{10}{\scriptsize $\Sigma_3$}

	\endrelabelbox
	
	\caption{\label{shrinkfig2} \small In the figure above, $\E_1$ has a convex shrinkwrapping surface $\Sigma_1$, and $\E_3$ has a concave shrinkwrapping surface $\Sigma_3$. $\E_2$ is geometrically finite, and $\Sigma_2=S_2$.}
	
\end{figure}


Now, let $L=\Sigma\cap\alpha_2$ be the coincidence set. Consider $r$-tubular neighborhood $N_r(L)$ of $L$ in $M$ for some small $r>0$. Let $T=\Sigma\cap N_r(L)$.
Let $T^+$ and $T^-$ be the components of $T-L$. Notice that $L\subset \overline{T^\pm}$. By Lemma \ref{thinlem}, $\Sigma$ has well-defined normal vector everywhere except $\Sigma\cap\alpha_2$. For $q\in \overline{T}^\pm$ define the normal vectors $N^\pm(q)$ to $\overline{T}^\pm$ at $q$, pointing towards to the unbounded component of $\E-\Sigma$. By \cite[Lemma 1.27]{CG}, $T$ has a well-defined tangent cone. For any $p\in L$, let $\alpha(p)$ be the dihedral angle between tangent planes of $T^+$ and $T^-$ at $p$, i.e. the complementary angle between $N^+(p)$ and $N^-(p)$. In particular, $\alpha(p)$ is the angle of the tangent cone of a great bigon in \cite[Lemma 1.27]{CG}. 

We claim that the dihedral angle $\alpha(p)\leq \pi$ for any $p\in L$. Assume that there exists a point $p_0 \in L$ with $\alpha(p_0)>\pi$. Then there exists a small neighborhood $I'$ of $p_0$ in $L$ such that for any point $p'\in I'$, $\alpha(p')>\pi$. Notice that for $r>0$ small, $T$ is a topological disk. Furthermore, $T$ is minimal except along $L$. Notice that $\Sigma$ is least area surface satisfying shrinkwrapping properties by Lemma \ref{shrinklem} - Condition (4). This implies $T$ is a least area disk with thin obstacle $L$. In particular, there cannot be a smaller area disk with the same boundary staying in the same side of $\alpha_2$. However, by construction of $T$, $L$ is a folding curve for this "least area" disk. We cannot push $T$ to convex side to decrease area along $L$ in general because of the shrinkwrapping condition. However, near $I'\subset L$, $\alpha_2$ stays in the concave side of $T$, so we can push $T$ along the folding curve $I'$ to the convex side without breaking the shrinkwrapping condition. Then, by \cite{MY}, we can decrease the area of $T$ by pushing $T$ along $I'$ to the convex side, and  by keeping the boundary of $T$ fixed. However, this contradicts to $T$ being least area satisfying the shrinkwrapping conditions. Hence, we conclude that $\alpha(p)\leq \pi$ along $L$, and hence the unbounded component of $\E-\Sigma$ is mean convex.

The similar idea works when $\Sigma$ is convex shrinkwrapping surface ($\Sigma\cap\alpha_1\neq\emptyset$) by changing the orientation of the normal vectors $N^\pm(p)$ towards to the bounded component of $\E-\Sigma$. Then, if there is an interval $I'\subset L$ where $\alpha(p)>\pi$ along $I'$, we can push $T$ to the other side of $\alpha_1$ and decrease the area. This gives a contradiction.  Again, we conclude that $\alpha(p)\leq \pi$ along $L$, and hence the bounded component of $\E-\Sigma$ is mean convex. Step 2 follows. \hfill $\Box$

\vspace{.2cm}

We repeat this process in every geometrically infinite end $\E_i$, and obtain a convex or concave shrinkwrapping surface $\Sigma_i$ in $\E_i$.  Recall that if we have a geometrically finite end $\E_j$, we take $\Sigma_j$ as convex surface which is in the boundary of the convex core in that end. 

In the trivial case (Case 3), $\Sigma$ does not touch the geodesics, then $\Sigma$ is already smoothly embedded closed minimal surface, and we are done. So, assuming the trivial case never happens, we have a convex or concave shrinkwrapping surface $\Sigma_i$ in every end $\E_i$ of $M$. 
The next two steps in the proof treat the two possible cases when every $\Sigma_i$ is concave,
or when at least one $\Sigma_i$ is convex.

\vspace{.2cm}

\noindent {\bf Step 3:} If $M$ contains at least one geometrically finite end, or one of the shrinkwrapping surfaces is convex, then $M$ contains a smoothly embedded minimal surface. 

\vspace{.2cm}

\noindent {\em Proof of Step 3:} Our aim is to cut the manifold $M$ from each end and cusp, and get a compact mean convex manifold $\wh{M}$. If we already have a convex surface in an end, we use that surface. If we do not have such convex surface in an end, or in a cusp, we choose a surface sufficiently far from the compact core, and change the metric near the boundary. Then, using a nontrivial isotopy or homology class of the convex surface in an end, we get an area minimizing surface $\Sigma$ in $\wh{M}$. Then, we show that $\Sigma$ is away from the regions where the metric was modified by using our bounded diameter lemmas. This shows $\Sigma$ is a smoothly embedded closed minimal surface in $M$.

Let $\C_M$ be the compact core of $M$. If $\E_i$ is a geometrically finite end, let $\wh{S}_i$ be the convex surface in the boundary of the convex core.

Similarly, if $\E_i$ is a geometrically infinite end which contains a convex shrinkwrapping surface $\Sigma_i$, then we take $\wh{S}_i=\Sigma_i$. 

If $\E_i$ contains a concave shrinkwrapping surface, then we choose $\wh{S}_i$ sufficiently far from the compact core $\C_M$ as follows. If $\E_i$ have bounded geometry, then let $C_i$ be the constant given by the bounded diameter lemma (Lemma \ref{boundiamlem}). Let $d(\wh{S}_i(0),\C_M)>C_i$. If $\E_i$ has unbounded geometry, let $S_i\times\{0\}\subset \partial \C_M$ where $\E_i\simeq S_i\times [0,\infty)$. Then, let $C_i=C_{K_i}$ be the constant for $K_i= S_i\times\{0\}$ in the end $\E_i$ given by the pointwise bounded diameter lemma (Lemma \ref{boundiamlem3}) . Choose $\wh{S}_i=S_i\times \{t\}$ such that $d(\wh{S}_i, S_i\times\{0\})>C_i$. 

If $M$ has cusps $X_1, X_2, ... ,X_m$, we modify the manifold $M$ as follows. Let $X_j\simeq T^2\times [0,\infty)$. Let $C_j$ be the constant in Lemma \ref{cusplem} for $X_j$. Then, let $\wh{T}_j=T^2\times \{2C_j\}$. 

Now, we define a compact mean convex manifold $\wh{M}$ as follows. For any end $\E_i$, cut $M$ through $\wh{S}_i$, and throw away the unbounded component of $\E_i-\wh{S}_i$. For any cusp $X_j$, cut $M$ through $\wh{T}_j$ and throw away the unbounded component of $X_j-\wh{T}_j$. We get a compact manifold $M'$ such that $\partial M'=\bigcup_{i=1}^N \wh{S}_i\bigcup_{j=1}^m\wh{T}_j$. For an end $\E_i$ with concave shrinkwrapping surface, modify the metric near $\wh{S}_i$ so that $M'$ becomes mean convex near $\wh{S}_i$. Similarly, for any $j$, modify the metric near $\wh{T}_j$ so that $M'$ becomes mean convex near $\wh{T}_j$. Let $\wh{M}$ be the compact mean convex manifold with this modified metric near the boundary.

We finish the proof in 3 cases: $M$ contains either at least 3 ends, or exactly 2 ends, or only one end.

\vspace{.2cm}

\noindent {\em Case 1: $M$ has at least 3 ends.} Recall that we have at least one  convex surface $\wh{S}_{i_0}$ in an end $\E_{i_0}$. Let $\xi$ be the homology class of $\wh{S}_{i_0}$. If $M$ has more than 2 ends, then $\xi$ is nontrivial in $H_2(\wh{M})$. Then, by Lemma \ref{LAlem}, we have a least area surface $\Sigma$ in the homotopy class of $S_{i_0}$. $\Sigma$ is nonempty as the homology class $\xi$ is nontrivial. Furthermore, as $M$ has more than 2 ends, for any surface $S'$ in the homotopy class of $S_{i_0}$, we have $S'\cap \C_M\neq \emptyset$. Hence, $\Sigma\cap \C_M \neq \emptyset$.

As $\wh{M}$ is mean-convex, $\Sigma$ is smoothly embedded in $\wh{M}$. Furthermore, since $\Sigma\cap \C_M \neq \emptyset$, then $d(\Sigma, \wh{S}_i)>C_i$ for any end $\E_i$ with concave shrinkwrapping surface. Similarly, as $\Sigma$ is area minimizing, $d(\Sigma, \wh{T}_j)>C_j$ by Lemma \ref{cusplem}. For the ends with convex surface, the boundary is already mean convex in the original metric. This implies $\Sigma$ is completely in the part of $\wh{M}$ where the metric is unmodified. Therefore, $\Sigma$ is a smoothly embedded minimal surface in the original manifold $M$. Hence, If $M$ has at least 3 ends, we are done.

In the following cases, we use a similar idea. All we need to make sure that the homotopy class of $S_{i_0}$ is not trivial, and the least area representative $\Sigma$ intersect the compact core $\C_M$.

\vspace{.2cm}

\noindent {\em Case 2: $M$ has exactly 2 ends.} Assume that $M$ is not a product manifold $S\times \BR$ (Type I Exceptional Manifold). Then, $M$ has a non-degenerate compact core $\C_M$. Here, nondegenerate means $\C_M$ must be $3$-dimensional, codimension-$0$ compact set in $M$. Then, we can employ the same idea as in the previous step, by cutting the manifold from noncompact parts in the ends, and cusps (if exists), and get a mean convex compact $3$-manifold $\wh{M}$. The homology class $\xi$ of the convex shrinkwrapping surface $\wh{S}_{i_0}$ is not trivial as $M$ has two ends. Then, the area minimizing surface $\Sigma$ in the homology class of $\xi$ in $\wh{M}$ will be a smoothly embedded closed minimal surface in the original manifold $M$ as before.


When $M$ is a product manifold $S\times \BR$ (Type I Exceptional Manifold), then the compact core $\C_M\simeq S\times\{t\}$ for any $t\in\BR$. Then, our method fails in the existence of concave shrinkwrapping surface in the other end. For details of this situation, see Section \ref{sec-Exc}.

\vspace{.2cm}

\noindent {\em Case 3: $M$ has exactly one end.} In this case, we have only one end $\E$ and it contains a convex surface $S$. Again, we cut our end $\E$ from $S$, and cut the cusp regions (if exists) sufficiently deep as before. We obtain a compact mean convex manifold $\wh{M}$. If $M$ has cusps, then the homology class of $S$ is nontrivial as before, and we get area minimizing surface $\Sigma$ in  $M$.

When $M$ is cusp-free, we need additional assumption that $M$ is not a handlebody (Type II Exceptional Manifold). As $M$ is not a handlebody, the isotopy class of $S$ is not trivial in $M$. In particular, if $\X_S$ is the space of surfaces in $\wh{M}$ isotopic to $S$, then $\inf\{|T| \mid T\in \X_S\}=c_0>0$. Then, by Lemma \ref{LAlem}, there exists a smoothly embedded least area surface $\Sigma$ in the compact mean convex manifold $\wh{M}$ in the isotopy class of $S$ achieving the infimum $c_0>0$, i.e. $|\Sigma|=c_0$. 

In each 3 cases, we obtained smoothly embedded closed minimal surface $\Sigma$ in $M$. Step 3 follows. \hfill $\Box$

\vspace{.2cm}
	
\noindent {\bf Step 4:} If every end $\E_i$ of $M$ contains a concave shrinkwrapping surface $\Sigma_i$ for $1\leq i \leq N$, then $M$ contains a smoothly embedded minimal surface. 

\vspace{.2cm}

\noindent {\em Proof of Step 4:} In this step, we use minmax lemma (Lemma \ref{minmaxlem}) to obtain a smoothly embedded minimal surface in $M$. In particular, we construct an open bounded set $\Omega$ in $M$ where $\overline \Omega$ is a strictly mean concave manifold. 


Now, by assumption, every end $\E_i$ of $M$ contains a concave shrinkwrapping surface $\Sigma_i$. In particular, they are minimal everywhere except some curve segments. Along these curve segments, the dihedral angle is less than $\pi$, i.e. convex towards the component containing infinity. Hence, when we evolve $\Sigma_i$ under the level set flow, by the main theorem (and Section 11) of \cite{Wh}, $\Sigma_i$ immediately becomes a smooth, strictly mean concave (normal vector points to $\infty$ side of $\E_i$) surface $\Sigma_i'$ in the positive side ($\infty$ side of $\E_i$) of $\Sigma_i$. We will take $\Sigma_i'$ as our smooth, strictly mean concave surfaces in the end $\E_i$ for each $i$. If there are cusps $X_j=T^\times[0,\infty)$, let $T_j=T^2\times\{C_j\}$ for sufficiently large $C_j$. Notice that $T_j$ has constant mean curvature $1$, and hence mean concave. Let $Z=\bigcup_i\Sigma_i'\bigcup_j T_j$.

Now, $Z$ separates a bounded open domain $\Omega$ in $M$. Furthermore, $\overline{\Omega}$ has smooth strictly concave boundary $Z$. Hence, by  Lemma \ref{minmaxlem} and Corollary \ref{minmaxcor}, $M$ contains a smoothly embedded unstable minimal surface $\Sigma$ with $\Sigma\cap\Omega\neq \emptyset$. Step 4 follows. \hfill $\Box$.


This proves $M$ contains a smoothly embedded closed minimal surface unless it is  an exceptional manifold. The proof of the theorem follows.
\end{pf}

\begin{rmk} [Stability of the Minimal Surfaces] Notice that all the minimal surfaces constructed above are obtained by using area minimization except one case: Hyperbolic $3$-manifold with geometrically infinite ends with bounded geometry where all ends contains a concave shrinkwrapping surface (Step 4). In this case, we used a min-max method \cite{Mo,So}. Hence, this is the only case we get an {\em unstable minimal surface}. In all other cases, the minimal surfaces we obtain are {\em stable}. 
\end{rmk}

\section{Final Remarks} 

\subsection{Exceptional Manifolds} \label {sec-Exc} \

\vspace{.2cm}

In this part, we focus on the exceptional manifolds, and see what can be derived more about these cases by using the proof of the main result.

\vspace{.2cm}

\noindent $\diamond$  {\bf Type I: Geometrically Infinite Product Manifolds:} 

Let $M$ be hyperbolic $3$-manifold which is a product, i.e. $M\simeq S\times \BR$ for some closed surface of genus $\geq 2$. Then, in the following cases, we can deduce the existence of a smoothly embedded closed minimal surface in $M$ by employing the techniques in the proof of main theorem. 
\begin{itemize}
	\item $M$ is geometrically finite
	\item Both shrinkwrapping surfaces are convex.
	\item Both shrinkwrapping surfaces are concave.	
\end{itemize}

However, only in the case where one end contains a convex shrinkwrapping surface, and the other contains concave shrinkwrapping surface, our methods fail because of the following. The compact core $\C_M$ will be degenerate, and any $S\times \{t\}$ can be taken as the compact core. In this case, if one end has convex shrinkwrapping surface, and the other end has concave shrinkwrapping surface, then we cannot use the bounded diameter idea in Step 3 of Theorem \ref{mainthm}. This is because, the location of surface we cut the manifold from the end containing concave shrinkwrapping surface depends on the location of compact core. Since in the other cases, the area minimizing surface must intersect the "well-defined" codimension-$0$ compact core, we can use our bounded diameter lemmas to cut deep enough in the end where $\Sigma$ cannot reach that far. However, in the product case, we do not have a "well-defined" place for the compact core, and the idea fails. 

\vspace{.2cm}

\noindent $\diamond$  {\bf Type II: Handlebodies (Schottky Manifolds).} 

This case is quite interesting. First notice that as $M$ is a handlebody, $\pi_1(M)$ is a free group. Since we do not allow rank-$1$ cusps, all elements of $\pi_1(M)$ represented by hyperbolic isometries. So, $M$ must be a Schottky manifold \cite{Mas}.

Furthermore, in the trivial case $M=\BH^3$, we know that there is no smoothly embedded closed minimal surface in $M$. The reason for this, the concentric $r$-spheres $S_r=\partial\B_r(0)$ with mean curvature $\coth{r}>1$ foliate $\BH^3$. If $\Sigma$ was a closed, smoothly embedded minimal surface in $\BH^3$, then the largest $r$ with $S_r\cap \Sigma\neq \emptyset$ would give the first point of touch, and give a contradiction with maximum principle.

In another trivial case, we can also show $M$ contains no closed minimal surface. Geometrically simplest Schottky Manifolds are Fuchsian Schottky Manifolds. Because of their simple geometric structure, they admit a natural CMC foliation, which gives nonexistence by maximum principle.

\begin{lem} \label{FuchsianSchottky} Let $M$ be a Fuchsian Schottky $3$-manifold. Then, there exists no closed embedded minimal surface in $M$.
\end{lem}

\begin{pf} Since $M$ is Schottky, $\pi_1(M)\simeq \mathbf{F}_g$ where $\mathbf{F}_g$ is a free group with $g$ generators for some $g\geq 2$. In particular, $M$ is a genus $g$ handlebody, which is in the homotopy type of a disk with $g$ punctures.

As $M$ is also Fuchsian, then the representation of the fundamental group of the manifold is in $Isom(\BH^2)=PSL(2,\BR)$, i.e. $\pi_1(M) \hookrightarrow PSL(2,\BR)\subset PSL(2,\mathbb{C})$. Hence, there exists a hyperbolic plane $\p_0$ in $\BH^3$ fixed by $\Gamma$. Then, we have a completely geodesic surface $\Sigma_0=\p_0/\Gamma$ in $M$. Furthermore, let $\Lambda=\PI \p_0$ be the asymptotic boundary of $\p_0$ which is a round circle in $S^2_{\infty}(\mathbf{H}^3)$. Now, consider the foliation of $\BH^3$ by constant mean curvature $H$-planes $\{\p_H\}$ with $\PI \p_H=\Lambda$ where $H\in(-1,1)$ (Figure \ref{fig_schottky}-left).

\begin{figure}[t]
	\begin{center}
		$\begin{array}{c@{\hspace{.1in}}c@{\hspace{.1in}}c}

		\relabelbox  {\epsfxsize=1.5in \epsfbox{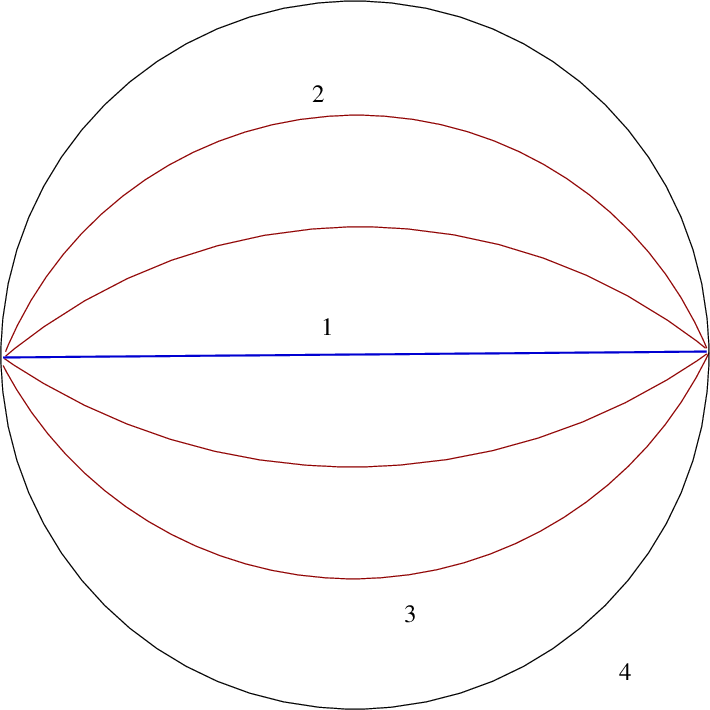}}
		\relabel{1}{\tiny $\p_0$} 
		\relabel{2}{\tiny $\p_H$}  
		\relabel{3}{\tiny $\p_{\mbox{-}H}$} 
		\relabel{4}{\tiny $\BH^3$} 
		\endrelabelbox &
		
		\relabelbox  {\epsfxsize=1.5in \epsfbox{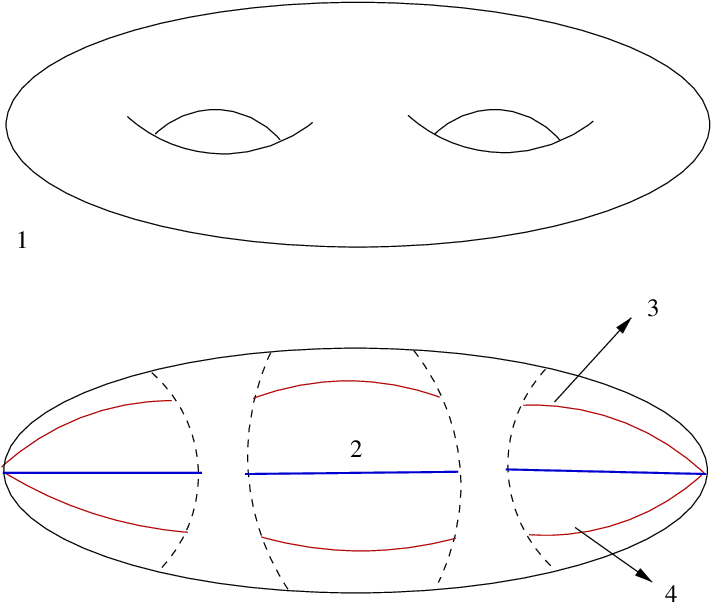}}
		\relabel{1}{\tiny $M$} 
		\relabel{2}{\tiny $\Sigma_0$}  
		\relabel{3}{\tiny $\Sigma_H$} 	
		\relabel{4}{\tiny $\Sigma_{\mbox{-}H}$} 	
		\endrelabelbox &
		
		\relabelbox  {\epsfxsize=1.5in \epsfbox{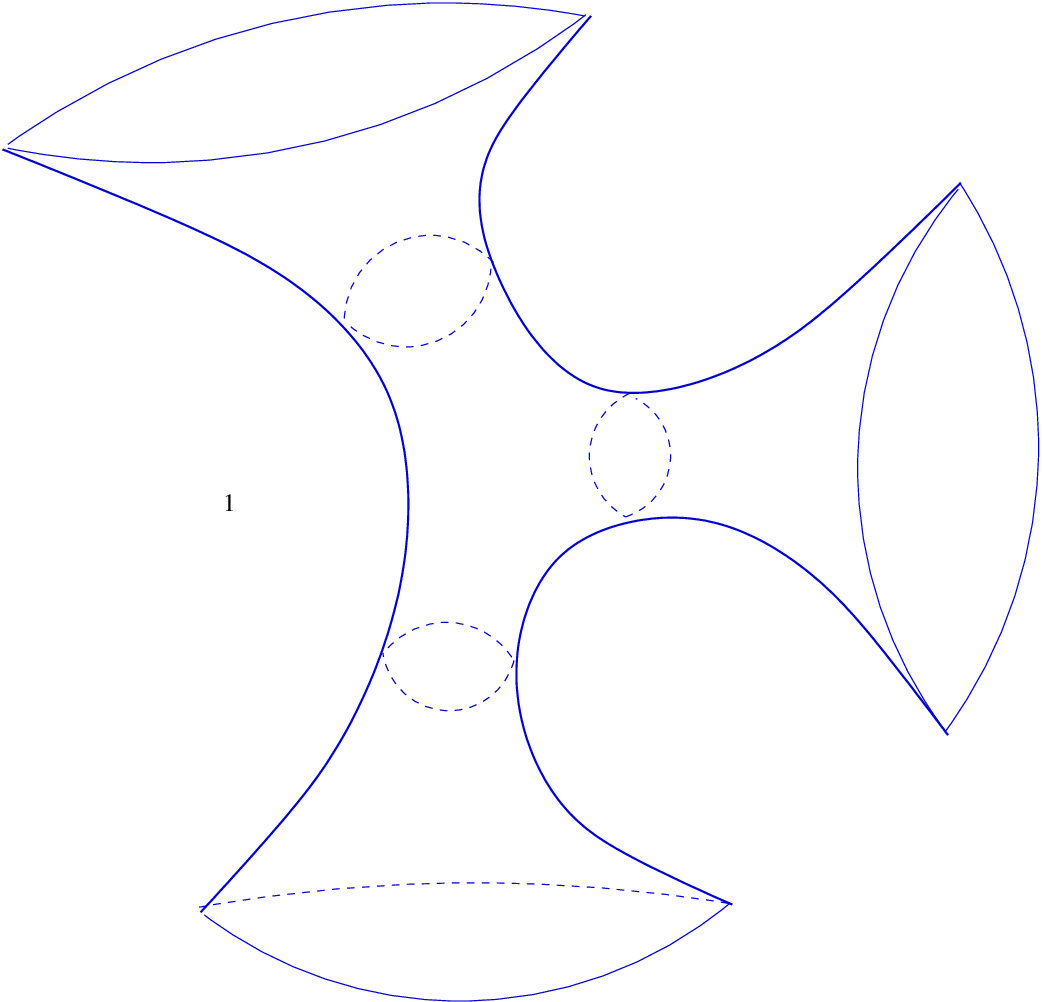}}
		\relabel{1}{\scriptsize $\Sigma_0$} 
		
		\endrelabelbox \\
		
	\end{array}$
	
\end{center}

\caption{ \label{fig_schottky} \scriptsize In the figure left, we have foliation of $\BH^3$ by $H$-planes for $H\in (-1,1)$. In the middle top, we have a topological picture of a Schotky manifold $M$, which is a handlebody. In the middle bottom, we have a picture of $M$ from the side. When we project the foliation $\{\p_H\}$ the manifold, we have a CMC foliation of $M$, i.e. $\pi(\p_H)=\Sigma_H$. In the right, we see the real geometric picture of the geodesic surface $\Sigma_0$ (hyperbolic 3-punctured sphere).}
\end{figure}

By construction, $\Lambda$ is preserved by $\Gamma$. Furthermore, as the manifold is oriented, the isometries cannot switch regions in the domain of discontinuity $\S^2_\infty(\mathbf{H}^3)-\Lambda=\Omega^+\cup\Omega^-$. $\p_H$ is the unique $H$-plane with $\PI\p_H=\Lambda$ for $H\in(-1,1)$. Then, $\p_H$ is preserved by $\Gamma$, i.e. for any $\varphi\in\Gamma$, we have $\varphi(\p_H)=\p_H$. This induces a foliation $\{\Sigma_H\mid H\in (-1,1)\}$ where $\Sigma_H=\p_H/\Gamma$. Intuitively, one can imagine this foliation as follows: $M$ is a genus $g$ handlebody, and  $\Sigma_0$ is the middle slice (largest cross-section) of $M$. Hence, $\Sigma_0$ is topologically a disk with $g$ holes. Then, upper half of $M-\Sigma_0$ is foliated by $\{\Sigma_H\mid H>0\}$, and the lower half of $M-\Sigma_0$ is foliated by $\{\Sigma_H\mid H<0\}$. See Figure \ref{fig_schottky}-middle.

Let $S$ be a closed embedded minimal surface in $M$. Consider the last point of touch $H_0=\sup\{H\mid \Sigma_H\cap S\neq \emptyset\}$. As $S$ is closed, $H_0<1$. If $H_0<0$, then we can use the last point of touch in the other side. So, we assume $H_0\in [0,1)$. However, this implies $\Sigma_{H_0}$ and $S$ has tangential intersection while one lies in the other side of the other. This is a contradiction by maximum principle. The proof follows.
\end{pf}

This proves that in the simplest class of Schottky manifolds, there is no closed embedded minimal surface. On the other hand, for a general Schottky manifold, while constructing a least area surface is not possible because of the trivial topology, it might still be possible to use some variation of min-max methods to construct a closed embedded minimal surface. However, as the ends are geometrically finite in Schottky Manifolds, the current min-max methods are not sufficient to obtain such a minimal surface.

\vspace{-.3cm}

\subsection{The Existence Problem and the Remanining Cases} \ \label{sec-remaining}

	\vspace{.1cm}
	
In this paper, we adressed the following question: {\em Which hyperbolic $3$-manifolds contain a closed embedded minimal surface?} After our results, there are only 3 cases left to finish off this problem: 

	\vspace{.1cm}
	
\noindent $\diamond$ {\bf Type I Exceptional Manifolds:} {\em Does every Type I exceptional manifold contain a closed embedded minimal surface?} One way to approach to this question is to use Gromov's dichotomy. Recently, Antoine Song proved a dichotomy suggested by Gromov \cite{So, Gr}: A complete, non-compact $3$-manifold either contains a closed, embedded minimal surface, or it contains a mean convex foliation. Hence, one needs to rule out existence of mean convex foliation in this case. One can visualize this situation as follows. $M=S\times \BR$ where the diameter of $|S\times \{t\}|= D(t)$ for $D:\BR\to \BR^+$ is a monotone decreasing function. The area minimizing limit escapes to infinity, and level sets are arranged so that they are mean convex. In such a case, $M$ cannot contain a closed minimal surface. 

On the other hand, in this question, the problem can be split into two cases. One case is that $M$ has bounded geometry, and the other is that $M$ has unbounded geometry. If one wants to construct closed embedded minimal surface in $M$, these two cases might need completely different techniques.
	
	\vspace{.1cm}
	
\noindent $\diamond$ {\bf  Type II Exceptional Manifolds:} {\em Does there exist a Type II exceptional manifold which contains a closed embedded minimal surface?} Our example of Fuchsian Schootky manifolds, and $\BH^3$ are the only nonexistence examples so far. However, it might be possible to construct similar CMC foliations with noncompact surfaces in such manifolds in general. Then, this would imply that no Type II exceptional manifold contains a closed embedded minimal surface.

	\vspace{.2cm}
	
\noindent $\diamond$ {\bf Rank-1 Cusps:} While rank-$2$ cusps have finite volume, and their structure is of the form $\T^2\times[0,\infty)$, rank-$1$ cusps have infinite volume, and their structure is of the form $S^1\times \BR\times[0,\infty)$. In the existence of rank-$1$ cusps, our methods do not apply, as we do not have any form of cusp lemma which explains the situation of closed embedded minimal surfaces in these cusps. One important case to start with is the manifold $M\simeq S\times \BR$ where $S$ is a punctured surface. Again, there is no existence or nonexistence results even in this basic case for this type of manifolds.

\subsection{Open Embedded Minimal Surfaces.} \

\vspace{.2cm}

One can generalize the question at hand by allowing  complete, noncompact, properly embedded minimal surfaces. Hence, the following question becomes interesting: {\em Which hyperbolic 3-manifolds admit an open, embedded minimal surface?}

For this version, the ends of the minimal surface lives in the cusps, and in the ends of the hyperbolic $3$-manifold. For cusped case, one can obtain embedded, finite area, noncompact minimal surfaces \cite{CHMR}. For the convex cocompact (geometrically finite and cusp-free) hyperbolic $3$-manifolds, \cite{AM} gives a positive answer to this question in various cases.


\begin{thebibliography}{MSY}

\bibitem[Ag]{Ag} I. Agol, Tameness of hyperbolic 3–manifolds, eprint math.GT/0405568.


\bibitem[AM]{AM} S. Alexakis,   and R. Mazzeo, {\em
Renormalized area and properly embedded minimal surfaces in hyperbolic 3-manifolds}, Comm. Math. Phys. {\bf 297} (2010) 621--651. 

\bibitem[An]{An} M. Anderson, {\em Complete minimal varieties in hyperbolic space}, Invent. Math. {\bf 69} (1982) 477--494.

\bibitem[Bo]{Bo} F. Bonahon, {\em Bouts des vari'et'es hyperboliques de dimension 3 (Ends of hyperbolic 3-manifolds)}, Ann. of Math. {\bf 124} (1986) 71--158.

\bibitem[BO]{BO} F. Bonahon, and J.-P. Otal, {\em Variétés hyperboliques à géodèsiques arbitrairement courtes}, Bull. London Math. Soc. {\bf 20}  (1988) 255--261.

\bibitem[BD]{BD} J.F. Brock, and N.M. Dunfield, {\em Norms on the cohomology of hyperbolic 3-manifolds}, Invent. Math. {\bf 210} (2017) 531--558.

\bibitem[BP]{BP} R. Benedetti, and C. Petronio, {\em Lectures on hyperbolic geometry}, Universitext. Springer-Verlag, Berlin, 1992. 

\bibitem[CG]{CG} D. Calegari, and D. Gabai, {\em Shrinkwrapping and the taming of hyperbolic 3-manifolds}, J. Amer. Math. Soc. {\bf 19} (2006) 385--446.

\bibitem[CL]{CL} G.R. Chambers, and Y. Liokumovich, {\em Existence of minimal hypersurfaces in complete manifolds of finite volume},  Invent. Math. \textbf{219} (2020) 179--217.

\bibitem[CT]{CT} J.W. Cannon, and W.P. Thurston, {\em Group invariant Peano curves}, Geom. Topol. {\bf 11} (2007) 1315--1355. 

\bibitem[CHMR1]{CHMR} P. Collin, L. Hauswirth, L. Mazet and H. Rosenberg, {\em Minimal surfaces in finite volume hyperbolic $3$-manifold}, Trans. Amer. Math. Soc. {\bf 369} (2017) 4293--4309.


\bibitem[CHMR2]{CHMR2} P. Collin, L. Hauswirth, L. Mazet and H. Rosenberg, {\em Erratum to ”Minimal surfaces in finite volume noncompact hyperbolic 3-manifolds}, arXiv:1810.07562


\bibitem[Co1]{Co1} B. Coskunuzer, {\em Embedded Plateau Problem}, Trans. Amer. Math. Soc. 364 (2012) 1211-1224.

\bibitem[Co2]{Co2} B. Coskunuzer, {\em Embedded $H$-Planes in Hyperbolic 3-Space}, Trans. Amer. Math. Soc. \textbf{371} (2019) 1253–1269.



\bibitem[Fe]{Fe} H. Federer, {\em Geometric measure theory}, Springer-Verlag, New York 1969.

\bibitem[FHS]{FHS} M. Freedman, J. Hass, and P. Scott, {\em Least area incompressible surfaces in 3-manifolds}, Invent. Math. {\bf 71} (1983) 609--642.



\bibitem[Gr]{Gr} M. Gromov. {\em Plateau-Stein manifolds}, Cent. Eur. J. Math. \textbf{12} (2014) 923--951.

\bibitem[Ha1]{Ha} J. Hass, {\em Minimal Fibrations of Hyperbolic $3$-manifolds}, arXiv:1512.04145.

\bibitem[Ha2]{Ha2} J. Hass, {\em Surfaces minimizing area in their homology class and group actions on 3-manifolds}, Math. Z. {\bf 199} (1988) 501--509.

\bibitem[HS]{HS} J. Hass, and P. Scott, {\em The existence of least area surfaces in 3-manifolds}, Trans. Amer. Math. Soc. {\bf 310} (1988) 87--114.
 
\bibitem[HW1]{HW} Z. Huang, and B. Wang, {\em Minimal fibration in hyperbolic $3$-manifolds fibering over the circle}, Proc. London Math. Soc. {\bf 118} (2019) 1305--1327.

\bibitem[HW2]{HW2} Z. Huang, and B. Wang, {\em Closed minimal surfaces in cusped hyperbolic $3$-manifolds}, Geom. Dedicata (2017) {\bf 189} 17--37.


\bibitem[IMN]{IMN} K. Irie, F. C. Marques, and A. Neves, {\em Density of minimal hypersurfaces for generic metrics}, Ann. of Math. {\bf 187} (2018) 963--972.

\bibitem[La]{La} M. Lackenby, {\em Heegaard splittings, the virtually Haken conjecture and property ($\tau$)}, Invent. Math. {\bf 164} (2006) 317--359. 


\bibitem[Mar]{Ma} A. Marden, {\em Hyperbolic manifolds. An introduction in 2 and 3 dimensions}, Cambridge University Press, Cambridge, 2016.

\bibitem[Mas]{Mas} B. Maskit, {\em Kleinian groups}, Springer-Verlag, Berlin, 1988. 

\bibitem[Mi]{Mi} Y.N. Minsky, {\em End invariants and the classification of hyp. $3$-manifolds}, Current developments in mathematics, 2002, 181--217, Int. Press, Somerville, MA, 2003. 



\bibitem[Mo]{Mo} R. Montezuma, {\em Min-max minimal hypersurfaces in non-compact manifolds}, J. Diff. Geom. {\bf 103} (2016) 475--519.

\bibitem[MR]{MR}  L. Mazet and H. Rosenberg, {\em Minimal hypersurfaces of least area}, J. Differential Geom. {\bf 106} (2017) 283--316.

\bibitem[MSY]{MSY} W.H. Meeks, L. Simon, and S.T. Yau, {\em Embedded minimal surfaces, exotic spheres, and manifolds with positive Ricci curvature}, Ann. of Math. {\bf 116} (1982) 621--659.

\bibitem[MY]{MY} W.H. Meeks, and S.T. Yau, {\em The classical Plateau problem and the topology of three-dimensional manifolds}, Topology {\bf 21} (1982) 409--442.

\bibitem[Pi]{Pi} J. Pitts, {\em Existence and regularity of minimal surfaces on Riemannian manifolds}, Mathematical Notes 27, Princeton University Press, Princeton (1981).


\bibitem[So1]{So1} A. Song,  {\em Existence of infinitely many minimal hypersurfaces in closed manifolds}, arXiv:1806.08816

\bibitem[So2]{So} A. Song, {\em A dichotomy for minimal hypersurfaces in manifolds thick at infinity}, arXiv:1902.06767



\bibitem[Su]{Su} D. Sullivan, {\em A finiteness theorem for cusps}, Acta Math. {\bf 147} (1981) 289--299.

\bibitem[Th]{Th} W.P. Thurston, {\em The geometry and topology of 3-manifolds}, Princeton University Lecture Notes, online at http://www.msri.org/publications/books/gt3m, 1982.


\bibitem[Wh]{Wh} B. White, {\em The size of the singular set in mean curvature flow of mean-convex sets}, J. Amer. Math. Soc. {\bf 13} (2000) 665--695.


\end{thebibliography}
\end{document}